\documentclass[12pt]{article}
\usepackage{authblk}
\author[1]{D. Gintides}
\author[2]{L. Mindrinos}
\affil[1]{\textit{\small Department of Mathematics, National Technical University of Athens, Greece, dgindi@math.ntua.gr} }
\affil[2]{\textit{\small Department of Natural Resources Development and Agricultural Engineering, Agricultural University of Athens, Greece, leonidas.mindrinos@aua.gr}}

\title{On the inverse elastic problem for isotropic media using Eshelby and Lippmann-Schwinger integral formulations}
\date{}

\usepackage{amsfonts}
\usepackage{graphicx}
\usepackage{amsmath}

\usepackage{amssymb,bm}
\usepackage{color}

\usepackage[colorlinks = true,
            linkcolor = blue,
            urlcolor  = blue,
            citecolor = blue,
            anchorcolor = blue]{hyperref}

\DeclareMathAlphabet\mathbit
    \encodingdefault\rmdefault\bfdefault\itdefault
\DeclareOldFontCommand{\bi}{\normalfont\bfseries\itshape}{\mathbit}

\newcommand{\be}{\begin{equation}}
\newcommand{\ee}{\end{equation}}

\def\fakebold#1{\relax\ifvmode\leavevmode\fi%
\ifmmode%
\setbox0=\hbox{$#1$}%
\else%
\setbox0=\hbox{#1}%
\fi%
\kern-.02em\copy0 \kern-\wd0%
\kern .04em\copy0 \kern-\wd0%
\kern-.0125em\raise.02em\box0%
}%

\renewcommand{\geq}{\geqslant}
\renewcommand{\leq}{\leqslant}
\newcommand{\R}{\mathbb{R}}
\newcommand{\C}{\bm{\mathcal{C}}}
\newcommand{\W}{\bm{\mathcal{W}}}
\renewcommand{\L}{\bm{\mathcal{L}}}
\newcommand{\V}{\bm{\mathcal{V}}}
\newcommand{\di}{\bm{\hat{d}}}

\newcommand{\n}{\bm{\hat{n}}}
\newcommand{\xhat}{\hat{\bm{x}}}

\usepackage{amsfonts}
\usepackage{graphicx}
\usepackage{amsmath}
\usepackage{amssymb} 
\newtheorem{theorem}{Theorem}

\usepackage{caption}
\usepackage{subcaption}

\begin{document}

\maketitle 
\begin{abstract}

We present two applications of the integro-differential volume equation for the eigenstrain, building on Eshelby's inclusion method \cite{eshelby1957, eshelby1961}, in the contexts of both static and dynamic linear elasticity. The primary objective is to address the inverse problem of recovering the elastic moduli of the inhomogeneity using a limited number of incident fields. In the static case, we adopt an efficient reformulation of Eshelby's equation proposed by Bonnet \cite{bonnet17}. By employing a first-order approximation in addition with a limited number of incident loadings and measurements, we numerically determine the material coefficients of the inclusion. In elastodynamics, we focus on the inverse scattering problem, utilizing the Lippmann-Schwinger integral equation to reconstruct the elastic properties of the inclusion through a Newton-type iterative scheme. We construct the Fr\'echet derivative and we formulate the linearized far-field equation. Additionally, the corresponding plane strain problems are analyzed in both static and dynamic elasticity.
\end{abstract}

 
 \section{Introduction}
 
 Inverse scattering problems in linear elasticity are of significant importance across various research domains. From a theoretical perspective, they pose complex questions for mathematicians, such as ensuring the uniqueness and stability of solutions, as well as addressing challenging numerical and computational problems, including the treatment of singular integrals. From a practical point of view, elasticity serves as a fundamental concept in nearly all branches of engineering, with applications in fields such as seismology, materials science, and the imaging of biological tissues.

Inverse problems are broadly classified into two categories: the inverse obstacle problem, which seeks to determine the geometry (position and shape) of a homogeneous elastic medium from measured data and the inverse medium problem, which focuses on reconstructing the material parameters of an inhomogeneous medium based on available data \cite{bonnet05, ColKre13}. Both categories refer to ill-posed problems and measured data can either be near-field data, such as Cauchy data collected on the accessible boundary of the medium, or far-field data, represented by the far-field patterns of the scattered field.

In the context of mathematical modeling, the Lam\'e equation serves as the fundamental framework, formulated within a bounded domain together with appropriate boundary conditions. Numerous variations arise depending on the specific setup, such as whether the problem is three- or two-dimensional, static or dynamic, and on the elastic properties of the medium, if it is homogeneous or inhomogeneous, isotropic or anisotropic. The choice of boundary conditions corresponds to different types of media, such as rigid scatterers, cavities, or inclusions.

In this work, we investigate the inverse medium problem in linear elasticity, focusing on the reconstruction of material parameters from far-field data. The primary focus is on the inverse scattering problem in $\R^3$ for an anisotropic elastic medium in the context of elastodynamics. Additionally, we address several simplified cases, including the static problem, isotropic media, and the corresponding plane strain problem, to provide a more comprehensive analysis.

Concerning unique solvability, no theoretical results exist without imposing conditions on the Lam\'e parameters, and most of the them rely on the knowledge of all boundary data. For foundational work on this topic, we refer to the series of papers by Nakamura and Uhlmann \cite{nakamura93, nakamura94, nakamura95}. On the other hand, numerous studies address the parameter identification problem, primarily through the use of linearization schemes. For isotropic media, regularization methods have been applied using internal data in both the dynamic  \cite{Bal} and the static case \cite{Hubmer}, as well as using boundary data for the time-dependent  \cite{armin15} and the static problem \cite{MarKouCha24}. Direct reconstruction formulas for 3D isotropic media in elastodynamics, under the assumption of known density, are presented in \cite{Bar07}. For anisotropic media, integral equation methods are employed for both the static case \cite{bonnet17, GinKir15} and the dynamic case \cite{bonnet2016}. Newton-type methods have also been explored for the reconstruction of small inclusions in $\R^2$  \cite{Feh06} and for experimental data in layered elastic media \cite{Doy00}. These methods often require the computation of the Fr\'echet derivative and its adjoint for the parameter-to-solution mapping, as discussed in \cite{ikehata, Kirsch16}.

In \cite{GinKir15}, the volume integral equation formulation was applied to the isotropic case, while in \cite{bonnet17}, it was extended to account for anisotropic elastic properties. This approach utilized Eshelby's equivalent inclusion method \cite{eshelby1957, eshelby1961}. For static elasticity, an integro-differential equation was derived, and its solvability was analyzed. The core idea of this method is to construct an equivalent problem in which the inhomogeneity is replaced by an inclusion that shares the same material properties as the surrounding exterior domain.

In \autoref{sec_static_all}, we revisit the results from \cite{bonnet17} concerning the anisotropic elastostatic case in three dimensions. Following this approach, we also address the corresponding plane strain problem, presenting the necessary adjustments and the different forms of the singular integrals. The elastodynamic problem is examined in \autoref{sec_dynamic_all} for both two- and three-dimensional setups. The direct scattering problem is reformulated as a Lippmann-Schwinger-type integral equation, and its well-posedness in Hilbert spaces is established. Later, we tackle the inverse medium problem, aiming to recover the Lam\'e parameters and the mass density as functions of the spatial variable. Specifically, we focus on reconstructing the material parameters from the far-field patterns of scattered waves generated by a small number of incident waves. The Fr\'echet derivative of the forward mapping is defined and characterized, and the linearized far-field integral equation is derived and solved using classical Tikhonov regularization.

In the final section, we present few numerical examples to demonstrate the applicability of the proposed scheme. These examples focus on the static case in two dimensions for an isotropic medium. We focus on the reconstruction of one of the Lam\'e parameters ($\lambda$) assuming a Gaussian-type representation. This part of our work should be seen as an initial step towards addressing the significantly more complex full reconstruction problem in three dimensions for anisotropic media.

\section{The static elastic inhomogeneity problem }\label{sec_static_all}

\subsection{Three-dimensional setup}\label{sec_3d_static}

Let $\Omega \subset \R^3 $ denote a bounded domain occupying the inhomogeneity.  Assume that $\R^{3} 
\setminus \overline{\Omega}$ is filled up with an isotropic and homogeneous elastic medium specified by the Lam\'{e} parameters $\lambda ,\, \mu,$ where $3\lambda
+2\mu >0$ and $\mu >0$.  The material elastic fourth-order tensor in $\mathbb{R}^{3} \setminus
\overline{\Omega}$ is denoted by $\C=\left( \C_{ijkl}\right) $ where 
\begin{equation*}
\C_{ijkl}=\lambda ~\delta _{ij}\delta _{kl}+\mu \left( \delta _{ik}\delta
_{jl}+\delta _{il}\delta _{jk}\right) ,~~~i,j,k,l=1,2,3,
\end{equation*}
satisfying the relations \cite{mura87}
\begin{equation}\label{eq_symmetry}
\C_{ijkl}=\C_{jikl}=\C_{ijlk}=\C_{lkij}.
\end{equation}

Assume that $\Omega$ is a bounded, simply connected domain in $\mathbb{R}^{3}$, with a smooth $C^{2}$ and closed boundary $\partial \Omega,$ having  different elastic properties from the outer medium, denoted by $\C^\ast$ satisfying also \eqref{eq_symmetry} and the strong ellipticity condition:
\begin{equation*}
\bm\xi : \C^\ast (\bm x) : \bm\xi > c_0 |\bm\xi |^2 >0,
\end{equation*}
for almost all $\bm x \in\Omega$ and all symmetric $3\times 3$ real matrices $\bm\xi,$ and $c_0>0.$ Here $``:"$ denotes the double inner product. 

We define the material contrast
\begin{equation}\label{eq_contrast}
\Delta \C (\bm x) := \C^\ast (\bm x) -\C,
\end{equation}
for $\mu^\ast \neq \mu$ and $\kappa^\ast \neq \kappa,$ where $\kappa$ is the bulk modulus and the strain tensor
\[
\bm\epsilon (\bm u) := \frac{1}{2} \left( \nabla \bm u + \nabla \bm u^\top  \right),
\]
for a vector-valued field $\bm u.$

Let $\bm u^{inc} (\bm x), \, \bm x \in \R^3$ denote the applied displacement field, not assumed to satisfy the Navier equation in free space. Then by $\bm u^{+} (\bm x), \, \bm x \in \R^3 \setminus \overline{\Omega}$ we define the total field in the exterior domain and by $\bm u^{-} (\bm x), \, \bm x \in \Omega$ the perturbation field.

We consider the following exterior transmission boundary value
problem: Find  $\bm u^+ \in C^2(\mathbb{R}^{3} \setminus \overline{\Omega})\cap C^1(\mathbb{R}^{3} \setminus {\Omega})$ in the interior domain $\Omega$  and 
 $\bm u^-  \in C^2(\Omega) \cap C^1( \overline{\Omega})$ in the
 exterior domain $\mathbb{R}^{3} \setminus \overline{\Omega}$ satisfying:
\begin{subequations}\label{eq_system}
\begin{alignat}{3}
\nabla\cdot \left( \C :\bm\epsilon (\bm u^+)  \right) &= \nabla\cdot \left( \C :\bm\epsilon (\bm u^{inc})  \right), \quad && \mbox{in }  \R^3 \setminus \overline{\Omega}, \label{eq_system1}\\
\nabla\cdot \left( \C^\ast :\bm\epsilon (\bm u^-)  \right) &= \nabla\cdot \left( \C^\ast :\bm\epsilon (\bm u^{inc})  \right), \quad && \mbox{in }  \Omega, \label{eq_system2}\\
\bm u^+ &= \bm u^- , \quad && \mbox{on }  \partial\Omega, \\
\bm T (\bm u^+) &= \bm T^\ast (\bm u^-) , \quad && \mbox{on }  \partial\Omega, \\
\bm u^+ - \bm u^{inc} &= \mathcal{O}\left(\frac{1}{|\bm x|} \right), \quad && \mbox{as } |\bm x| \rightarrow \infty, \label{eq_system5}
\end{alignat}
\end{subequations}
where the traction operators are defined by
\begin{align*}
\bm T (\bm u^+) &= \left[\C : \bm\epsilon (\bm u) |_{+} \right] \cdot \n (\bm x) , \\
\bm T^\ast (\bm u^-) &= \left[\C^\ast : \bm\epsilon (\bm u) |_{-} \right] \cdot \n (\bm x),
\end{align*}
for the unit normal vector $\n$ to $\partial\Omega$ pointing in $\R^3 \setminus \overline{\Omega}.$ Here $(\cdot)|_{\pm}$ denotes the limiting values on $\partial\Omega.$ In addition, the mechanical strain can be derived via an eigenstrain which can be produced e.g.  by  a thermal cause and it is not derived by a displacement field \cite{mura87}. In this case, we have to modify equations \eqref{eq_system1} --\eqref{eq_system2}.

The asymptotic behavior \eqref{eq_system5} holds uniformly over all angles in $\mathbb{R}^{3} \setminus \overline{\Omega}$ is necessary for the well-posedness of the above problem \cite{Kup79}.  

For the application of the Eshelby's method \cite{eshelby1957, eshelby1961} the fundamental tensor (Green's function) $\bm G$ of the outer space is needed. Its analytic form is given by
\cite[Eq. 5.8]{mura87} 
\begin{equation}
\bm G (\bm x) = \lambda ^{\prime }\frac{1}{\left\vert
\bm x\right\vert } \bm I +\mu ^{\prime } \hat{\bm{x}} \otimes \hat{\bm{x}}, \quad \hat{\bm{x}} = \frac{\bm x}{|\bm x|},
\end{equation}
where  $\bm I$ is the second-order identity and 
\begin{equation}\label{eq_l_m}
\lambda ^{\prime }=\frac{\lambda +3\mu }{4\pi \mu \left( \lambda +2\mu
\right) }, \quad \mu ^{\prime }=\frac{\lambda +\mu }{4\pi \mu \left( \lambda
+2\mu \right) }. 
\end{equation}
Here the symbol $\otimes $ denotes juxtaposition of two vectors, that is the  well-known dyadic product of two vectors in $\mathbb{R}^{n}$ and summation convention implied over the repeated indices.

We define the volume vector potential
\begin{equation}\label{eq_volume}
\begin{aligned}
\W (\bm h) (\bm x) &:= \nabla\cdot \int_{\R^3} \bm G (\bm x- \bm y) \cdot \bm h (\bm y) d \bm y, \\
&=  \int_{\R^3} \nabla \bm G (\bm x- \bm y) \cdot \bm h (\bm y) d \bm y,
\end{aligned}
\end{equation}
for a symmetric tensor density $\bm h.$ It is a bounded operator seen as a mapping $\W : L^2_c (\R^3;\R^{3\times 3}_{sym}) \rightarrow H^1_{loc} (\R^3;\R^3)$ and if $D:= \mbox{supp} (\bm h),$ then the displacement field $\bm w = \W (\bm h)$ satisfies \cite{bonnet17}
\begin{subequations}\label{eq_equi_system}
\begin{alignat}{3}
-\nabla\cdot \left( \C :\bm\epsilon (\bm w)  \right) &= 0, \quad && \mbox{in }  \R^3 \setminus \overline{D}, \\
-\nabla\cdot \left( \C :\bm\epsilon (\bm w)  \right) &=  \nabla\cdot \bm h, \quad && \mbox{in }  D, \label{eq_equi_system2}\\
\bm w^+ &= \bm w^- , \quad && \mbox{on }  \partial D, \\
\bm T (\bm w^+) &= \bm T (\bm w^-) + \bm h \cdot \n, \quad && \mbox{on }  \partial D, \\
\bm w (\bm x) &= \mathcal{O}\left(\frac{1}{|\bm x|} \right), \quad && \mbox{as } |\bm x| \rightarrow \infty.
\end{alignat}
\end{subequations}

From \eqref{eq_equi_system} we have that a field $\bm v$ satisfying the equation $\bm v = \bm u^{inc} + \W (\Delta \C : \bm\epsilon (\bm v))$ solves also \eqref{eq_system}. Thus, the transmission problem \eqref{eq_system} admits the equivalent form of a system of volume integro-differential equations of Lippmann-Schwinger type:
\begin{equation}\label{eq_sol}
\bm u^{\mp}  = \bm u^{inc} +\W (\Delta \C : \bm\epsilon (\bm u^-)), \quad \mbox{in } \Omega \cup (\R^3 \setminus \overline{\Omega}). 
\end{equation}
We first solve the volume  integro-differential equation for $\bm u^-$  in  $\Omega$ an then we use the integral representation for the field $\bm u^+$.

In addition to the  the Lippmann-Schwinger formulation we employ an Eshelby approach to solve the direct problem based on \cite{bonnet17}.  We define a simpler inclusion problem, which is well-posed as a classical potential problem with body forces for the displacement field $\bm w$ \cite{Kup79}.  We have to point out that now the material parameters of $D$ are those of the exterior domain $\mathbb{R}^{3} \setminus \overline{\Omega}$.  Then, we derive an equation for the interior strain field by applying constant invertible  tensors from the left.

Following \cite{bonnet17}, equation \eqref{eq_sol} can be reformulated as a Fredholm integral equation of the second kind  where the resulting integral operator is a contraction mapping for an equivalent strain, as discussed above.  The following theorem provides a summary of this procedure, presented without proof.

\begin{theorem}
The volume integral equation \eqref{eq_sol} for $\bm x \in \Omega$ is written equivalently  as
\begin{equation}\label{eq_sol_fred}
(\bm{\mathcal{I}} -\bm{\mathcal{Q}}) \bm h^\ast = 2 (\C^\ast_B  +\bm{\mathcal{I}})^{-1} : \bm{\mathcal{B}}^{-1} : \Delta \C : \bm\epsilon (\bm u^{inc}),
\end{equation}
for the modified unknown density function $\bm h^\ast := \bm{\mathcal{B}}^{-1} : \Delta \C : \bm\epsilon (\bm u^-),$ 
where
\begin{align*}
\bm{\mathcal{Q}} &= (\C^\ast_B  +\bm{\mathcal{I}})^{-1} : (\C^\ast_B  -\bm{\mathcal{I}}) : (\bm{\mathcal{I}} +2 \bm{\mathcal{B}}:\bm{\mathcal{H}}:\bm{\mathcal{B}}), \\
\C^\ast_B &= \bm{\mathcal{B}}^{-1} : \C^\ast : \bm{\mathcal{B}}^{-1}, \\
\bm{\mathcal{H}} &= \bm\epsilon (\W), \\
\bm{\mathcal{B}} &= \sqrt{2\mu} \left[ \tfrac{1}{3} \left( \sqrt{\tfrac{1+\nu}{1-2\nu}} -1  \right) \bm I \otimes \bm I + \bm{\mathcal{I}} \right]
\end{align*}
and $\nu$ is the Poisson's ration and $\bm{\mathcal{I}}$ is the fourth-order identity tensor. 

The operator $\bm{\mathcal{Q}} : L^2_c (B;\R^{3\times 3}_{sym}) \rightarrow L^2_c (B;\R^{3\times 3}_{sym})$ is a contraction ($\parallel \bm{\mathcal{Q}} \parallel<1$) and therefore $\bm{\mathcal{I}} - \bm{\mathcal{Q}}$ is invertible with bounded inverse. Thus, the exterior transmission problem \eqref{eq_system} admits a unique solution which can be constructed using Neumann series expansion. 

\end{theorem}

\subsection{Plane strain problem}

We now consider the modifications required to apply the same approach in \(\R^2\). Following \cite{GinKir15} and we focus on the differences due to the dimension reduction. First of all, the fundamental solution is  given by
\begin{equation}\label{green2d}
\bm\Phi (\bm x) = \lambda ^{\prime } \ln \left(\frac{1}{\left\vert
\bm x\right\vert } \right) \bm I +\mu ^{\prime } \hat{\bm{x}} \otimes \hat{\bm{x}}, \quad \hat{\bm{x}} = \frac{\bm x}{|\bm x|}, \quad \bm x \in \R^2,
\end{equation}
and the Lam\'e parameters are now required to satisfy $\mu >0$ and $\lambda + \mu \geq 0.$ 

The transition from \eqref{eq_system2} to \eqref{eq_equi_system2} is due to Eshelby's idea to define $\bm h$ such that:
\begin{equation}
 \C^\ast :\bm\epsilon (\bm u^- - \bm u^{inc})  =  \C :\bm\epsilon (\bm u^- - \bm u^{inc})  + \bm h, 
\end{equation}
or equivalently 
\begin{equation}
\bm h -  \Delta \C :\bm\epsilon (\bm u^- )  =  \Delta\C :\bm\epsilon ( \bm u^{inc}).
\end{equation}
The last equation using the definition \eqref{eq_volume} results in the following Fredholm integral equation of the second kind
\begin{equation}\label{eq_fred}
\bm h -  \Delta \C : \nabla\nabla \int_\Omega \bm\Phi (\bm x- \bm y)  : \bm h (\bm y) d \bm y =  \Delta\C :\bm\epsilon ( \bm u^{inc}),
\end{equation}
for a strongly singular integral operator, since from the analytic form of the Green function (see \eqref{green2d}) we have the asymptotic estimate   
\[
\left| \frac{\partial^2}{\partial x_i \partial x_j} \bm\Phi_{kl} (\bm x-\bm y) \right| \leq \frac{c}{|\bm x - \bm y|^2}, \quad i,j,k,l=1,2.
\]

Because of the singularity at $\bm x = \bm y$ we derive (for a sufficiently smooth $\bm h$) the formula  \cite[Eq. (20)]{niwa1986}
\begin{equation}\label{eq_jump}
 \nabla\nabla \int_\Omega \bm\Phi (\bm x- \bm y)  : \bm h (\bm y) d \bm y =  -\bm\gamma : \bm h +   \int_\Omega \nabla\nabla \bm\Phi (\bm x- \bm y)  : \bm h (\bm y) d \bm y,
\end{equation}
where
\[
\bm\gamma_{ijkl} = \frac{3\lambda + 7\mu}{8 \mu (\lambda+2\mu)} \delta _{ij}\delta _{kl} - \frac{\lambda + \mu}{8 \mu (\lambda+2\mu)}  \left( \delta _{ik}\delta
_{jl}+\delta _{il}\delta _{jk}\right) ,\quad i,j,k,l=1,2.
\]
We define 
\begin{equation}\label{eq_volume2}
\bm A (\bm h) (\bm x) :=  \int_{\Omega} \nabla\nabla \bm\Phi (\bm x- \bm y) : \bm h (\bm y) d \bm y
\end{equation}
and \eqref{eq_fred} using \eqref{eq_jump} admits the compact form
\begin{equation}\label{eq_eshelby1}
\left(\Delta\C^{-1} + \bm\gamma \right) : \bm h - \bm A (\bm h) = \bm\epsilon ( \bm u^{inc}).
\end{equation}

In the above equation we have used that $\Delta\C$ is invertible which is true given the definition \eqref{eq_contrast} and the assumptions on the parameters.

Next step is to show that the constant tensor $\bm M := \Delta\C^{-1} + \bm\gamma$ is invertible also in two dimensions. The elastic tensor is given by
\begin{equation}
\C = 2 \kappa \bm J + 2\mu \bm K, \quad \kappa = \frac{\mu}{1-\nu},
 \end{equation}
where
\[
\bm J_{ijkl} = \frac{1}{2} \delta _{ij}\delta _{kl}, \quad \bm{\mathcal{I}}_{ijkl} = \frac{1}{2}\left( \delta _{ik}\delta
_{jl}+\delta _{il}\delta _{jk}\right), \quad i,j,k,l=1,2
\]
and $\bm K = \bm{\mathcal{I}} -  \bm J.$ Then, the application of the Sherman-Morrison formula \cite{li2008} gives
\begin{equation}
\Delta\C^{-1} = \frac{1}{2 (\kappa^\ast - \kappa)} \bm J + \frac1{2(\mu^\ast - \mu)} \bm K. 
\end{equation}
The tensor $\bm{\mathcal{\gamma}}$ can also be rewritten  using the same decomposition so that invertibility of $\bm M$ is guaranteed. It is easy to see that
\[
\bm{\mathcal{\gamma}} = \frac{2\lambda + 6\mu}{4 \mu (\lambda+2\mu)} \bm J - \frac{\lambda + \mu}{4 \mu (\lambda+2\mu)} \bm K.
\]
Then,
\[
\bm M =  \alpha \bm J  + \beta \bm K,\]
where
\begin{equation}\label{eq_a_b}
\begin{aligned}
\alpha &= \frac{(\lambda+3\mu)(\kappa^\ast - \kappa)+\mu (\lambda+2\mu)}{2(\kappa^\ast - \kappa)(\lambda + 2\mu)},\\
\beta &= \frac{2\mu (\lambda+2\mu) - (\mu^\ast -\mu)(\lambda + \mu)}{4(\mu^\ast - \mu)(\lambda + 2\mu)}.
\end{aligned}
\end{equation}
Since $\mu >0, \, \lambda + \mu \geq 0$ and $\kappa^\ast \neq \kappa, \, \mu^\ast \neq
 \mu,$ we see that both denominators of $\alpha$ and $\beta$ are non-zero. We address the following cases:
 \begin{itemize}
 \item If $\kappa^\ast > \kappa \, \Rightarrow \, \alpha >0.$
 \item If $\kappa^\ast < \kappa,$ then  $\alpha \neq 0 \, \Leftrightarrow \, \kappa - \kappa^\ast \neq \dfrac{\mu (\lambda + 2\mu)}{\lambda + 3\mu}.$
 \item If $\mu^\ast < \mu \, \Rightarrow \, \beta <0.$
 \item If $\mu^\ast > \mu,$ then  $\beta  \neq 0 \, \Leftrightarrow \, \mu^\ast - \mu \neq \dfrac{2\mu (\lambda + 2\mu)}{\lambda + \mu}.$
 \item If $\lambda + \mu =0 \, \Rightarrow \, \beta \neq 0.$
 
 \end{itemize}

 Thus, only if $\kappa - \kappa^\ast = \tfrac{\mu (\lambda + 2\mu)}{\lambda + 3\mu}$ and $\mu^\ast - \mu = \tfrac{2\mu (\lambda + 2\mu)}{\lambda + \mu},$ we get $\alpha = \beta
=0$ and $\bm M$ will be the zero tensor transforming \eqref{eq_volume2} into a Fredholm integral equation of the first kind.
If $\alpha \neq 0$ or $\beta \neq 0,$ then the inverse exist. If both of them are non-zero then the inverse admits the form
\[
\bm M^{-1} = (\Delta\C^{-1} + \bm\gamma)^{-1} =  \frac1{\alpha} \bm J  + \frac1{\beta} \bm K.
\]

The kernel of the singular integral operator \eqref{eq_volume2} can be represented using its symmetrized version \cite{corna}
\begin{equation}\label{eq_symm}
\bm S_{ijkl} =\frac{1}{4} \left( \bm D_{ik,jl} + \bm D_{il,jk} +\bm D_{jk,il} +\bm D_{jl,ik}  \right),\quad i,j,k,l, =1,2,
\end{equation}
where
\[
\bm D_{ik,jl} = \frac{\partial^2}{\partial x_j \partial x_l} \bm\Phi_{ik}, 
\]
for the Green tensor (see \eqref{green2d})
\[
\bm\Phi_{ik} = \lambda' \ln \left(\frac{1}{|\bm x - \bm y|} \right) \delta_{ik} + \mu' \frac{(\bm x_i - \bm y_i)(\bm x_k - \bm y_k)}{|\bm x - \bm y|^2}.
\]
The last step is to show that the operator $\bm M-\bm A$ is Fredholm with index zero. To do so, we have to compute the determinant of its symbol tensor, defined by $\bm\Psi,$ and demonstrate that it is strictly positive (under certain assumptions on the material parameters). Then, the application of Noether's theory will justify the well-posedness of \eqref{eq_eshelby1} as in the three-dimensional case \cite{GinKir15}.

The symbol is defined by \cite{mikhlin87}
\begin{equation}
\bm\Psi (\bm k) = \bm M - \hat{\bm S}(\bm k),
\end{equation}
where $ \hat{\bm S} := \mathcal{F}\{\bm S\}$ denotes the  2D Fourier transform
  \[
  \mathcal{F}\{f\}(\bm k) =   \int_{\R^2}  f(\bm x) e^{i \bm k \cdot \bm x }d \bm x, \quad \bm k\in\R^2.
  \]
 We express the symbol as a $3\times 3$ matrix using Mandel notation
    \cite{bli96}:
    \begin{equation}
[\bm\Psi (\bm k)]_{3\times 3} = \begin{pmatrix}
\bm\Psi_{1111} & \bm\Psi_{1122} & \sqrt{2}\bm\Psi_{1112} \\
\bm\Psi_{1122} & \bm\Psi_{2222} & \sqrt{2}\bm\Psi_{2212} \\
\sqrt{2}\bm\Psi_{1112} & \sqrt{2}\bm\Psi_{2212} & 2\bm\Psi_{1212} 
\end{pmatrix},
\end{equation}
where
\[
\bm\Psi_{ijkl}(\bm k) = \bm M_{ijkl} - \hat{\bm S}_{ijkl}(\bm k), \quad i,j,k,l=1,2,
\]
for
\begin{equation}
\begin{aligned}
\bm M_{ijkl} &= \frac{(\alpha-\beta)}{2}\delta _{ij}\delta _{kl}   +\frac{\beta}2 \left( \delta _{ik}\delta
_{jl}+\delta _{il}\delta _{jk} \right), \\
\hat{\bm S}_{ijkl} &=-\frac{1}{4} \left( \bm k_j \bm k_l \hat{\bm\Phi}_{ik} + \bm k_j \bm k_k \hat{\bm\Phi}_{il} + \bm k_i \bm k_l \hat{\bm\Phi}_{jk} + \bm k_i \bm k_k \hat{\bm\Phi}_{jl}  \right).
\end{aligned}
\end{equation}
In the last equation, we have used that 
\begin{equation}
 \hat{\bm D}_{ik,jl} = \mathcal{F}\left\{\frac{\partial^2}{\partial \bm x_j \partial \bm x_l} \bm\Phi_{ik}\right\} = (-i \bm k_j)(-i \bm k_l) \hat{\bm\Phi}_{ik} =- \bm k_j \bm k_l \hat{\bm\Phi}_{ik}.
 \end{equation} 
 The functions $\hat{\bm\Phi}_{ik}, \, i,k=1,2$ are given by \cite{vladi}
  \begin{equation}
  \begin{aligned}
  \hat{\bm\Phi}_{11} (\bm k) &=   \lambda' \frac{1}{\bm k_1^2+ \bm k_2^2}  + \mu' \frac{\bm k_2^2-\bm k_1^2}{(\bm k_1^2+\bm k_2^2)^2}, \\
   \hat{\bm\Phi}_{12} (\bm k) &=   -2 \mu' \frac{\bm k_1 \bm k_2}{(\bm k_1^2+\bm k_2^2)^2}, \\
    \hat{\bm\Phi}_{22} (\bm k) &=   \lambda' \frac{1}{\bm k_1^2+\bm k_2^2}  + \mu' \frac{\bm k_1^2-\bm k_2^2}{(\bm k_1^2+\bm k_2^2)^2}.
  \end{aligned}
  \end{equation}

The determinant, after some lengthy calculations, admits the form
\begin{equation}
\det ([\bm\Psi (\bm k)]_{3\times 3}) = \alpha \beta^2 + \alpha \beta \lambda' + \frac{\beta^2}{2}(\lambda' - \mu') + \frac{\alpha+\beta}{4}(\lambda' - \mu')(\lambda' +\mu').
\end{equation}
As in the 3D case \cite{GinKir15} the determinant is independent of $\bm k$ and depends only on the physical parameters $\alpha, \,\beta$ defined in \eqref{eq_a_b} and on $\lambda' , \, \mu'$ given by \eqref{eq_l_m} and satisfying $\lambda'>0, \, \mu' \geq 0,$ and $\lambda'-\mu' = \tfrac{1}{2\pi (\lambda+2\mu)}>0.$ Thus, the sign of the determinant depends on the signs of $\alpha$ and $\beta.$ For instance, if $\kappa^\ast >\kappa$ and $\mu < \mu^\ast < 3\mu,$ then $\alpha,\, \beta>0$ and $\det ([\bm\Psi (\bm k)]_{3\times 3})>0.$ To summarize, we state that it holds $\inf |\det ([\bm\Psi (\bm k)]_{3\times 3})|>0$ almost always, meaning when $\bm M$ is invertible.

\section{The elastodynamic inhomogeneity problem}\label{sec_dynamic_all}

\subsection{Three-dimensional problem}\label{sec_3d_inhomo}

We consider the scattering problem of a time-harmonic incident field $\bm u^{inc},$ either $p$ or $s$ wave, by an elastic inhomogeneity occupying a bounded domain $\Omega \subset \R^3.$ The elastic parameters are defined as in \autoref{sec_3d_static}. Then, the scattered field  $\bm u^{sc}$ and the generated interior field $\bm u^{int}$ satisfy:
\begin{subequations}\label{eq_system_dynamic}
\begin{alignat}{3}
\nabla\cdot \left( \C :\bm\epsilon (\bm u^{sc})  \right) + \rho \omega^2 \bm u^{sc}  &= 0, \quad && \mbox{in }  \R^3 \setminus \overline{\Omega}, \label{eq_system1d}\\
\nabla\cdot \left( \C^\ast (\bm x) :\bm\epsilon (\bm u^{int})   \right)  + \rho^\ast (\bm x) \omega^2 \bm u^{int} &= 0, \quad && \mbox{in }  \Omega, \label{eq_system2d}\\
\bm u^{tot} &= \bm u^{int} , \quad && \mbox{on }  \partial\Omega, \\
\bm T (\bm u^{tot}) &= \bm T^\ast (\bm u^{int}) , \quad && \mbox{on }  \partial\Omega, \\
|\bm x | \left( \frac{\partial \bm u_{\alpha}^{sc}}{\partial |\bm x|} - i k_\alpha u_{\alpha}^{sc} \right) &\rightarrow 0,  \quad && \mbox{as } |\bm x| \rightarrow \infty,  \label{eq_system5d}
\end{alignat}
\end{subequations}
for $\alpha = p$ or $s.$ Here, $\omega$ is the angular frequency, $\rho^\ast >0$ and $\rho >0$ denote the mass density of $\Omega$ and  $\R^3 \setminus \overline{\Omega},$ respectively. The wavenumbers are given by
\begin{equation}\label{eq_wavenumbers}
k_p^2 = \frac{\rho \omega^2}{\lambda+2\mu}, \quad k_s^2 = \frac{\rho \omega^2}{\mu}.
\end{equation}
The radiation condition \eqref{eq_system5d} holds uniformly in all directions and it is valid since the exterior domain is isotropic. 

The problem \eqref{eq_system_dynamic}, the so-called direct scattering problem, is well-posed \cite{bonnet2016} and admits the equivalent integral formulation
\begin{subequations}\label{eq_ls_full}
\begin{alignat}{3}
\bm u^{tot} - \L (\bm u^{int}) &= \bm u^{inc}, \quad && \mbox{in }  \R^3 \setminus \overline{\Omega}, \label{eq_utot_3d}\\
\bm u^{int} - \L (\bm u^{int}) &= \bm u^{inc}, \quad && \mbox{in }  \Omega, \label{eq_LS_3d}
\end{alignat}
\end{subequations}
where
\[
\L (\bm v)(\bm x) := \W_\omega (\Delta \C : \bm\epsilon (\bm v)) (\bm x)+ \omega^2 \V_\omega (\Delta\rho \,\bm v)(\bm x), \quad \bm x \in \left(\R^3 \setminus \overline{\Omega}\right) \cup \Omega,
\]
with
\begin{equation}\label{eq_volume_dynamic}
\begin{aligned}
\V_\omega (\bm h) (\bm x) &:=  \int_{\Omega} \bm G_\omega (\bm x- \bm y) \cdot \bm h (\bm y) d \bm y, \\
\W_\omega (\bm h) (\bm x) &:= \nabla\cdot \int_{\Omega} \bm G_\omega (\bm x- \bm y) \cdot \bm h (\bm y) d \bm y.
\end{aligned}
\end{equation}
In the special case of $\rho^\ast =\rho$,  the basic integro-differential operator $\L$ simplifies to \[ \L (\bm v)(\bm x) = \W_\omega (\Delta \C : \bm\epsilon (\bm v)) (\bm x). \]

The fundamental tensor $\bm G_\omega$ for the dynamic case is given by \cite{Kup79}
\[
\bm G_\omega (\bm x) = -\frac{k_p}{\mu k_s^2} \nabla\nabla \frac{e^{i k_p |\bm x|}}{4\pi k_p |\bm x|} + \frac{1}{\mu k_s} \left( k_s^2 \bm I +\nabla\nabla \right) \frac{e^{i k_s |\bm x|}}{4\pi k_s |\bm x|}.
\]

The Lippmann-Schwinger type integral equation \eqref{eq_LS_3d} is uniquely solvable since the operator $\bm{\mathcal{I}} -\L$ seen as mapping from $H^1 (\Omega)$ into 
$H^1 (\Omega)$ is invertible with bounded inverse \cite{bonnet2016}. Then, the total field $\bm u^{tot}$ in $\R^3 \setminus \overline{\Omega}$ is given with respect to $\bm u^{int}$ by the integral equation \eqref{eq_utot_3d}. The analogous for isotropic inclusions follows from  \cite{GinKir15}.

Using the identity 
\[
\W_\omega (\bm h) (\bm x) = \V_\omega (\nabla\cdot\bm h) (\bm x) - \bm{\mathcal{S}}_\omega (\bm h \cdot \n) (\bm x),
\]
where
\[
\bm{\mathcal{S}}_\omega (\bm h ) (\bm x) = \int_{\partial\Omega} \bm G_\omega (\bm x- \bm y) \cdot \bm h (\bm y) ds (\bm y),
\]
is the well-known single-layer potential, we may rewrite the  operator $\L$ in the form
\begin{equation}\label{eq_op_L}
\L (\bm v)(\bm x) =   \V_\omega \left(\nabla\cdot (\Delta \C : \bm\epsilon (\bm v)) + \omega^2\Delta\rho \,\bm v\right)(\bm x) - \bm{\mathcal{S}}_\omega ((\Delta \C : \bm\epsilon (\bm v)) \cdot \n) (\bm x).
\end{equation}

For the corresponding inverse problem we will need the form  of the far-field patterns of the scattered field. Recall that any solution of \eqref{eq_system1d} satisfying \eqref{eq_system5d} admits the asymptotic behavior
\begin{equation}
\bm u^\infty (\bm x) =  \bm u^\infty_p (\hat{\bm{x}}) \frac{e^{i k_p |\bm x|}}{4\pi k_p |\bm x|} + \bm u^\infty_s (\hat{\bm{x}}) \frac{e^{i k_s |\bm x|}}{4\pi k_s |\bm x|} + \mathcal{O} \left(\frac{1}{|\bm x|}  \right),
\end{equation}
 as $|\bm x|\rightarrow \infty,$ uniformly in all directions $\hat{\bm{x}} = \bm x / |\bm x|\in S^2,$ where $S^2$ is the unit sphere. The pair $( \bm u^\infty_p, \,  \bm u^\infty_s)$ is called the far-field patterns.

The asymptotic behavior of the fundamental tensor together with \eqref{eq_op_L} result in the integral representation of the  far-field patterns \cite[Section 4.2]{bonnet2016}
\begin{equation}\label{eq_far}
\bm u^\infty_\alpha (\hat{\bm{x}}) = \bm F^\infty_\alpha (\hat{\bm{x}}), \quad \alpha=p,s,
\end{equation}
where
\begin{align*}
\bm F^\infty_p (\hat{\bm{x}}) &=  \frac{k_p}{\lambda+2\mu} \left(\bm I_p (\hat{\bm{x}}) \cdot \hat{\bm{x}} \right) \hat{\bm{x}}, \\
\bm F^\infty_s (\hat{\bm{x}}) &= \frac{k_s}{\mu} \left[\bm I_s (\hat{\bm{x}}) - \left( \bm I_s (\hat{\bm{x}}) \cdot \hat{\bm{x}} \right)  \hat{\bm{x}} \right], 
\end{align*}
with 
\begin{align*}
\bm I_\alpha (\hat{\bm{x}}) &= \int_\Omega e^{-i k_\alpha \hat{\bm{x}} \cdot \bm y} \left[  \nabla\cdot \left( \Delta \C : \bm\epsilon (\bm u^{int}) \right) +\omega^2 \Delta\rho  \,\bm u^{int} \right] d\bm x \\
&\phantom{=} - \int_{\partial\Omega} e^{-i k_\alpha \hat{\bm{x}} \cdot \bm y} \left( \Delta \C : \bm\epsilon (\bm u^{int}) \right) \cdot \n (\bm y) ds(\bm y), \quad \alpha = p,s.
\end{align*}

\subsubsection{The inverse problem}

We are interested in solving the inverse problem of reconstructing the material elastic tensor $\C^\ast,$ meaning the Lam\'e parameters $\lambda^\ast, \, \mu^\ast$ and the mass density $\rho^\ast$ in $\Omega$ from the knowledge of the far-field patterns $( \bm u^\infty_p (\hat{\bm{x}}), \,  \bm u^\infty_s (\hat{\bm{x}})),$ for all $\hat{\bm{x}} \in S^2,$ and few incident waves ($p$ or $s$). We vary the incident waves with respect to the direction of incidence $\hat{\bm{d}}_k, \,k=1,...,K$ and we denote them  by $\bm u_{\alpha,k}^{inc} (\bm x; \hat{\bm{d}}_k).$ In the following, we suppress the $k-$dependence for the sake of presentation.

The direct problem can be seen as the mapping 
\begin{equation}\label{eq_operator}
\bm{\mathcal{F}}: (L^\infty (\Omega))^3 \rightarrow (L^2 (S^2))^2, \quad 
 (\lambda^\ast, \, \mu^\ast, \, \rho^\ast) \mapsto ( \bm u^\infty_{p}, \,  \bm u^\infty_{s}).
\end{equation}

Thus, the corresponding inverse problem is to solve
\begin{equation}\label{eq_inverse}
\bm{\mathcal{F}} (\lambda^\ast, \, \mu^\ast, \, \rho^\ast) = ( \bm u^\infty_{p}, \,  \bm u^\infty_{s}) \quad \Rightarrow \quad (\bm F^\infty_p, \, \bm F^\infty_s) =  ( \bm u^\infty_{p}, \,  \bm u^\infty_{s}), \, \forall \hat{\bm{x}}\in S.
\end{equation}
 The operators $\bm F^\infty_\alpha$ depend non-linearly on the material parameters and thus we propose to linearize the above system of equations by taking the Fr\'echet derivative of $\bm{\mathcal{F}}.$ Next we show that the  Fr\'echet derivative, denoted by 
 \[
 \bm{\mathcal{F}}'(\lambda^\ast, \, \mu^\ast, \, \rho^\ast) (h_\lambda, \, h_\mu, \, q),
 \]
  is well defined and directly computed. The functions $(h_\lambda, \, h_\mu, \, q)$ are assumed to be sufficiently small.

From \eqref{eq_LS_3d} we see that the Fr\'echet derivative $\bm w := (\bm u^{int})'$ of $\bm u^{int}$ with respect to $(\lambda^\ast, \, \mu^\ast, \, \rho^\ast)$ in direction $(h_\lambda, \, h_\mu, \, q)$ satisfies the Lippmann-Schwinger integral equation
\begin{multline}\label{eq_LS_per}
\bm w (\bm x) -\W_\omega \left(\Delta \C : \bm\epsilon (\bm w)+\C_{h} : \bm\epsilon (\bm u^{int})\right) (\bm x) \\
- \omega^2 \V_\omega (\Delta\rho \,\bm w + q \, \bm u^{int} )(\bm x) = 0, \quad \bm x\in\Omega,
\end{multline}
where $\C_{h}$ denotes the perturbed elastic tensor. 

From  \eqref{eq_inverse} and \eqref{eq_LS_per} and the equivalence between the system of integral equations \eqref{eq_ls_full} and the systems of PDEs \eqref{eq_system1d}--\eqref{eq_system2d} the following theorem holds for the  Fr\'echet derivative of $\bm{\mathcal{F}}.$ 

\begin{theorem}
The operator \eqref{eq_operator} is Fr\'echet differentiable and its derivative is defined by
\[
 \bm{\mathcal{F}}'(\lambda^\ast, \, \mu^\ast, \, \rho^\ast) (h_\lambda, \, h_\mu, \, q) = ( \bm w^\infty_{p}, \,  \bm w^\infty_{s}),
\]
where $ ( \bm w^\infty_{p}, \,  \bm w^\infty_{s})$ denote the far-field patterns of the radiating solution $\bm w$ of
\begin{subequations}\label{eq_system_frechet}
\begin{alignat}{3}
\nabla\cdot \left( \C :\bm\epsilon (\bm w)  \right) + \rho \omega^2 \bm w  &= 0, \quad && \mbox{in }  \R^3 \setminus \overline{\Omega}, \label{eq_system1f}\\
\nabla\cdot \left( \C^\ast (\bm x) :\bm\epsilon (\bm w)   \right)  + \rho^\ast (\bm x) \omega^2 \bm w &= -\nabla\cdot \left( \C^\ast_h (\bm x) :\bm\epsilon (\bm u^{int})   \right) \nonumber\\
&\phantom{=} - q (\bm x) \omega^2 \bm u^{int}, \quad && \mbox{in }  \Omega. \label{eq_system2f}
\end{alignat}
\end{subequations}

\end{theorem}

Establishing the injectivity of the Fr\'echet derivative $\bm{\mathcal{F}}'$ seen as a mapping from parameters to far-field measurements is nontrivial,  despite various results available for forward operators that map parameters to boundary data. See, for example, \cite{nakamura93} for analogous results in the case $q=0$ with $h_\lambda$ and $h_\mu$ being smooth $C^\infty$ parameters. This task will be a subject of future work.

Let $\kappa := (\lambda^\ast, \, \mu^\ast, \, \rho^\ast)$ and $\bm{\mathcal{U}}^\infty := ( \bm u^\infty_{p}, \,  \bm u^\infty_{s}).$ We solve \eqref{eq_inverse} using a Newton-type method, meaning we consider the linearized equation
\begin{equation}\label{eq_reg1}
 \bm{\mathcal{F}}' (\tau_n) = \bm{\mathcal{U}}^\infty -  \bm{\mathcal{F}} (\kappa_n),
\end{equation}
to be solved for the update $\tau_n$ in order to obtain
$
\kappa_{n+1} = \kappa_n + \tau_n. $  
Here we have assumed that an initial guess is given. We handle the ill-posedness of \eqref{eq_inverse}  by a Tikhonov regularization scheme. Thus, the regularized solution of \eqref{eq_reg1} is taken as
\[
\kappa_{n+1} = \kappa_n + \left[\alpha_n \bm I + \bm{\mathcal{F}}'  (\kappa_n)^\ast \bm{\mathcal{F}}'  (\kappa_n) \right]^{-1} \bm{\mathcal{F}}'  (\kappa_n)^\ast ( \bm{\mathcal{U}}^\infty -  \bm{\mathcal{F}} (\kappa_n)),
\]
where $\alpha_n$ is the regularization parameter at the $n-$th step.

The convergence of the iterative scheme is guaranteed if the  non-linearity condition:
\[
\|\bm{\mathcal{F}} (\kappa) - \bm{\mathcal{F}} (\bar{\kappa}) - \bm{\mathcal{F}}' (\kappa) (\kappa - \bar{\kappa} ) \| \leq C \| \kappa - \bar{\kappa}\| \| \bm{\mathcal{F}} (\kappa) - \bm{\mathcal{F}} (\bar{\kappa}) \|,
\]
for all $\bar\kappa$ satisfying $\kappa=\bar\kappa,$ on $\partial\Omega,$ is satisfied. For a rigorous proof, in the case  of $q=0$ and boundary data, we refer to \cite{Hubmer} where the appropriate spaces are discussed. If $\bm{\mathcal{F}}$ models the so-called parameter-to-solution mapping (for all parameters and internal data), its Fr\'echet differentiability was examined in \cite{Kirsch16, armin15} for the time-dependent problem.

  \subsection{Plane strain problem}
  
  In this section we outline the differences of the two-dimensional problem compared to the three-dimensional one presented in \autoref{sec_3d_inhomo}. The time-harmonic incident field is either a longitudinal plane wave
  \[
\bm u^{inc}_p (\bm x ;  \di) = \di \,  e^{i k_p \di \cdot \bm x}, \quad \bm x \in\R^2,
\]
or a transversal plane wave
\[
\bm u^{inc}_s (\bm x ;  \di) = -\bm Q \cdot \di \,  e^{i k_s \di \cdot \bm x}, \quad \bm x \in\R^2,
\]
where $\di$ is the vector of incidence, the wavenumbers are given in \eqref{eq_wavenumbers} and  $\bm Q$ denotes the unitary matrix
\begin{equation*}
\bm Q= \begin{pmatrix}
\phantom{-}0 & 1\\ -1 & 0
\end{pmatrix}.
\end{equation*}

The system of equations \eqref{eq_ls_full} is valid also in this setup for the fundamental solution defined by
\begin{equation}\label{green2d_dynamic}
\begin{aligned}
\bm \Phi_\omega (\bm x) =  \frac{i}{4\mu}H_0^{(1)}(k_s \left|\bm x\right|)\bm I + \frac{i}{4\rho \omega^2} \nabla \nabla^\top \left[H_0^{(1)}(k_s \left|\bm x\right|)-H_0^{(1)}(k_p \left|\bm x \right|)\right],
\end{aligned}
\end{equation}
where $H_0^{(1)}$ denotes the  Hankel function  of order zero and of the first kind. Again, the asymptotic behaviour of the Hankel function results in the far-field representation of the integral operators \cite{ChaGinMin18, ChaKreMon00}.

We define the coefficients
\[
\beta_p = \frac{e^{i\pi/4}}{\lambda +2\mu}\frac1{\sqrt{8\pi k_p}}, \quad \beta_s = \frac{e^{i\pi/4}}{\mu}\frac1{\sqrt{8\pi k_s}}, 
\]
and the tensors  $\bm J_p (\xhat) = \hat{\bm{x}} \otimes \hat{\bm{x}}$ and $\bm J_s (\xhat) = \bm I - \bm J_p (\xhat).$  Then, the following representations hold
 \begin{equation*}
\begin{aligned}
\bm{\mathcal{S}}^\infty_{\omega,\alpha} (\bm h )(\xhat) &=\beta_\alpha \int_{\partial\Omega} \bm J_\alpha (\xhat) \cdot \bm h(\bm y) \, e^{-i k_\alpha \xhat \cdot \bm y}ds(\bm y), & \quad & \xhat \in S,\\
\V^\infty_{\omega,\alpha}( \bm h )(\xhat) &=\beta_\alpha \int_{\Omega} \bm J_\alpha (\xhat) \cdot \bm h(\bm y) \, e^{-i k_\alpha \xhat \cdot \bm y}d \bm y, & \quad & \xhat \in S
\end{aligned}
\end{equation*}
and as in \eqref{eq_far} we obtain
\begin{equation}\label{eq_far_2d}
\begin{aligned}
\bm u^\infty_\alpha (\hat{\bm{x}}) &= \V^\infty_{\omega,\alpha} (\nabla\cdot \left( \Delta \C : \bm\epsilon (\bm u^{int}) \right) +\omega^2 \Delta\rho  \,\bm u^{int} )(\xhat) \\
&\phantom{=} - \bm{\mathcal{S}}^\infty_{\omega,\alpha} (\left( \Delta \C : \bm\epsilon (\bm u^{int}) \right) \cdot \n (\bm y) )(\xhat), \quad \alpha=p,s,
\end{aligned}
\end{equation}
where $\bm u^{int}$ solves \eqref{eq_LS_3d}.

The terms in the integrands admit the forms
\begin{align*}
\Delta \C : \bm\epsilon (\bm u ) &= \Delta\lambda (\nabla\cdot \bm u ) \bm I + 2 \Delta\mu \,\bm\epsilon (\bm u ),  \\
\nabla\cdot \left( \Delta \C : \bm\epsilon (\bm u) \right) &=   \Delta\mu \Delta \bm u + (\Delta \lambda + \Delta\mu) \nabla (\nabla\cdot \bm u ) + (\nabla \Delta\lambda) (\nabla \cdot \bm u) \\
&\phantom{=}
  + 2 \bm\epsilon (\bm u ) \nabla \Delta \mu,
\end{align*}
for the differences of the Lam\'e parameters
\[
\Delta\lambda = \lambda^\ast (\bm x) - \lambda, \quad \Delta\mu = \mu^\ast (\bm x) - \mu.
\]

The definition of the direct and the corresponding inverse problem is exactly the same as in the 3D case and the linearization and regularization methods are formulated similarly.

\section{Numerical examples}

We examine the applicability of the proposed method in the two-dimensional regime. As a toy example we assume that $\mu^\ast = \mu$ and $\rho^\ast = \rho$ so that we focus on the reconstruction of the Lam\'e parameter $\lambda^\ast.$

\begin{figure}
\centering
\begin{subfigure}{.5\textwidth}
  \centering
  \includegraphics[width=.7\linewidth]{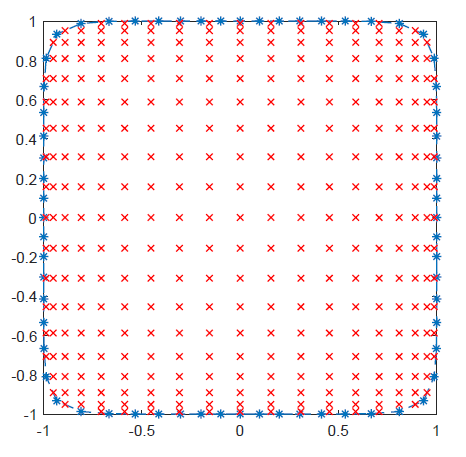}
\end{subfigure}%
\begin{subfigure}{.5\textwidth}
  \centering
  \includegraphics[width=.84\linewidth]{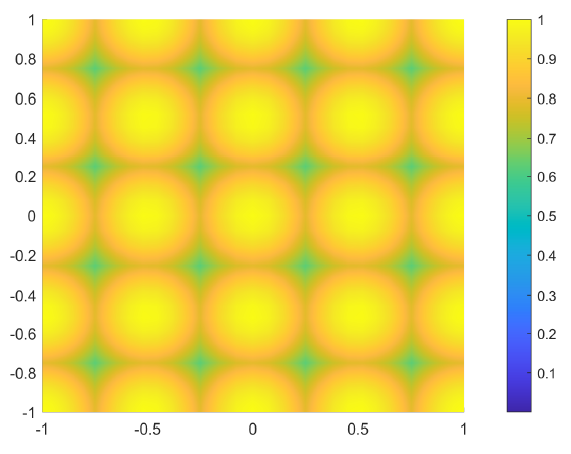}
\end{subfigure}
\caption{Left: The interior (red) and the boundary (blue) grid points for $n=32$ and $N=21$. Right: The Gaussian basis functions for $K=5.$}
\label{fig_setup}
\end{figure}

The computational domain is given by $D = [-1,1] \times [-1,1] \subset \R^2$ and $\Omega$ is a rounded square with boundary parametrization
\[
\partial\Omega = \{ \bm x (t) = q(t) (\cos t, \, \sin t ) : t \in [0,2\pi] \},
\]
for the radial function
$
q(t) = (\sin^{10} t + \cos^{10}t)^{-0.1}.
$
At the boundary we consider $2n$ equidistant points $t_j = \tfrac{j \pi}{n}, \, j=0,...,2n-1,$ whereas in $\Omega$ we use  non-uniformly distributed grid points 
\[
\bm x_{kj} := (x_1 (k),\,x_2 (j)) = -\left( \cos \left( \tfrac{k-1}{N-1}\pi\right),\,  \cos \left( \tfrac{j-1}{N-1}\pi\right)\right), \quad \text{for } k,j=1,...,N.
\]
In the left picture of Figure \ref{fig_setup} we present the boundary and domain nodal points for $n=32$ and $N=21.$ We approximate numerically the partial derivatives by the  differential quadrature method, meaning by a weighted sum at the nodal points \cite{BelKasCas72, BerMal96}. The volume integral in $\Omega$ is computed using the Clenshaw--Curtis quadrature method \cite{ClenCur60, Gen72a}.

We set $(\lambda,\,\mu)= (1,1)$ and $\rho =1.$ The far-field data $\bm{\mathcal{U}}^\infty (\hat{\bm{x}})$ are generated using \eqref{eq_far_2d} for four longitudinal incident waves of the form
  \[
\bm u^{inc}_{p,j} (\bm x ) = \di_j \,  e^{i k_p \di_j \cdot \bm x}, \quad j=1,2,3,4,
\]
for $\di_j = (\cos \bm\theta_j , \, \sin \bm\theta_j),$
where
$
\bm\theta = \left( 0,\, \tfrac{\pi}{2}, \,\pi, \, \tfrac{3\pi}{2}\right).
$
Then, we add noise with respect to the $L^2-$norm
\[
\bm{\mathcal{U}}^\infty_\delta = \bm{\mathcal{U}}^\infty + \delta \frac{\|\bm{\mathcal{U}}^\infty \|_2}{\| \bm{\mathcal{V}}\|_2}  \bm{\mathcal{V}}, 
\] 
where $\bm{\mathcal{V}}$ is a complex vector with components being normally distributed random variables and $\delta$ is a fixed noise level.

A Gaussian basis set is used to approximate the Lam\'e parameter $\lambda^\ast$ by
\[
\lambda^\ast (\bm x) \approx \sum_{k,j=1}^K \Lambda_{kj} \, e^{-\frac{K\pi}{4} \left[ (x_1 - z_k)^2 + (x_2 - z_j)^2\right]},
\]
where $z_k = -1 + \tfrac{2(k-1)}{K-1},$ for $k=0,1,...,K.$ Thus, the inverse problem is to recover the $K\times K$ coefficients $\Lambda. $ In the right picture of Figure \ref{fig_setup} we present the basis functions for $K=5.$

\begin{figure}
\centering
  \includegraphics[width=.8\linewidth]{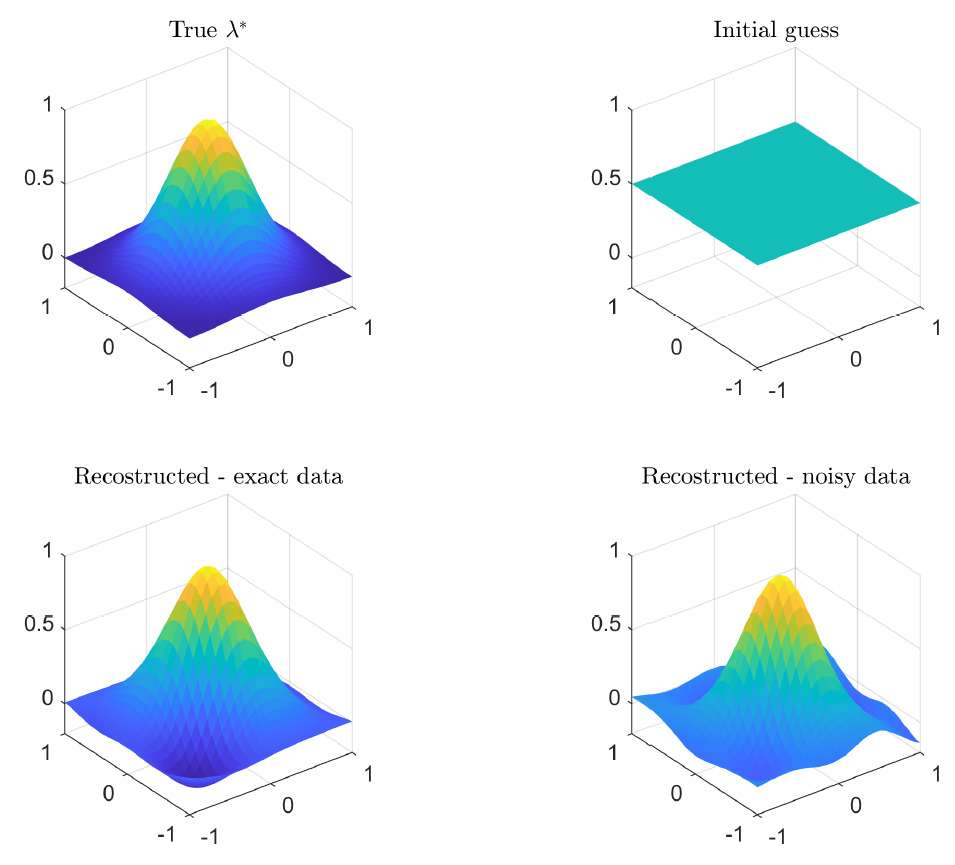}
\caption{The exact Lam\'e parameter $\lambda^\ast$ of the first example (top left), the initial guess $\lambda^\ast_0 = 0.5$ (top right) and the reconstructed values for exact (bottom left) and data with $3\%$ noise (bottom right).}
\label{fig_ex1}
\end{figure}

In the first example, the true function is given by
\[
\lambda^\ast_{ex} (\bm x) = e^{-\pi  |\bm x|^2}
\]
and we approximate it by $5\times 5$ Gaussian basis functions. We set $n=32, \, N=41$ for the number of nodal points and $\lambda=2$ is the parameter in the exterior domain and $\omega=0.1.$ Initially we assume a constant Lam\'e parameter $\lambda^\ast_0 = 0.5$ and the regularization parameter of the Tikhonov method is given by $\alpha_n = 0.001 \times 2^{-n},$ for $n=0,1,...$ at the $n-$th iteration step. The reconstructions after 2 iterations (exact data) and 4 iterations (noisy data with $\delta=3\%$) are presented in Figure \ref{fig_ex1}.

\begin{figure}
\centering
  \includegraphics[width=.8\linewidth]{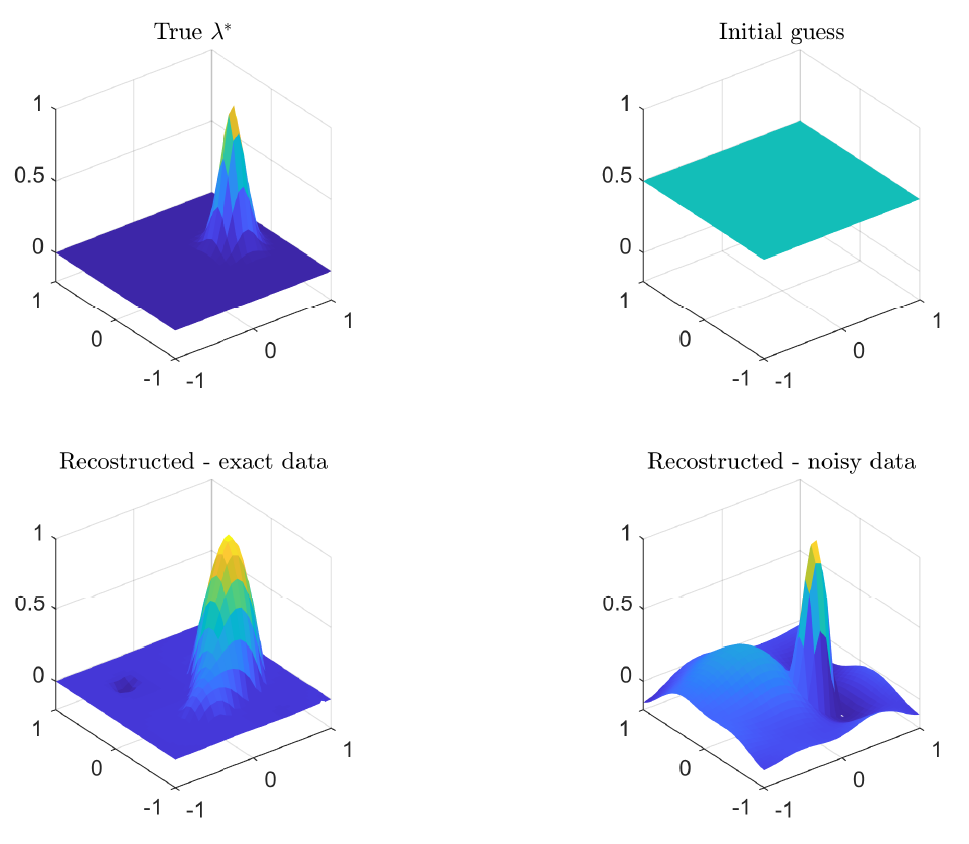}
\caption{The exact Lam\'e parameter $\lambda^\ast$ of the second example (top left), the initial guess $\lambda^\ast_0 = 0.5$ (top right) and the reconstructed values for exact (bottom left) and data with $3\%$ noise (bottom right).}
\label{fig_ex2}
\end{figure}

In the second example, we  aim to recover the function
\[
\lambda^\ast_{ex} (\bm x) = e^{- 15 \pi  |\bm x - (0.5,\,0)|^2}.
\]
We keep all parameters the same except for $\lambda=1.$
In Figure \ref{fig_ex2} we present the reconstructed parameters  after 3 iterations (exact data) and 4 iterations (data with noise).

In this section, we examined only a simplified version of the ``full" inverse problem to reconstruct both the Lam\'e parameters and the density function. This has to be seen as a proof of principle to justify the applicability of the proposed scheme. In future work we plan to investigate the full inverse problem in both two and three dimensions.
  
\section{Conclusions}

In this work, we studied inverse scattering problems in linear elasticity by addressing the reconstruction of material parameters for both isotropic and anisotropic media. In the static case, we applied the Eshelby's equivalent inclusion method to write the inverse problem as  an integro-differential equation for the eigenstrain. In elastodynamics, we reformulated the inverse scattering problem as a Lippmann-Schwinger-type integral equation.  The derivation and characterization of the Fr\'echet derivative of the forward operator, provided an efficient method for solving the linearized problem via Tikhonov regularization.  By using far-field patterns of scattered waves, we demonstrate the feasibility of recovering the material parameters, even with few incident waves. The study also includes numerical examples in simplified two-dimensional static scenarios for isotropic media, which serve as an initial step towards addressing the more complex full inversion in three-dimensional and anisotropic problems.


\bibliographystyle{siam}
\bibliography{ginmin_bib}

\begin{thebibliography}{10}

\bibitem{Bal}
{\sc G.~Bal, C.~Bellis, S.~Imperiale, and F.~Monard}, {\em Reconstruction of
  constitutive parameters in isotropic linear elasticity from noisy full-field
  measurements}, Inverse problems, 30 (2014), p.~125004.

\bibitem{Bar07}
{\sc P.~E. Barbone and A.~A. Oberai}, {\em Elastic modulus imaging: some exact
  solutions of the compressible elastography inverse problem}, Physics in
  Medicine \& Biology, 52 (2007), p.~1577.

\bibitem{BelKasCas72}
{\sc R.~Bellman, B.~G. Kashef, and J.~Casti}, {\em Differential quadrature: A
  technique for the rapid solution of nonlinear partial differential
  equations}, Journal of Computational Physics, 10 (1972), pp.~40--52.

\bibitem{BerMal96}
{\sc C.~W. Bert and M.~Malik}, {\em Differential quadrature method in
  computational mechanics: a review}, Applied Mechanics Reviews, 49 (1996),
  pp.~1--28.

\bibitem{bli96}
{\sc A.~Blinowski, J.~Ostrowska-Maciejewska, and J.~Rychlewski}, {\em
  Two-dimensional hooke's tensors-isotropic decomposition, effective symmetry
  criteria}, Archives of Mechanics, 48 (1996), pp.~325--345.

\bibitem{bonnet2016}
{\sc M.~Bonnet}, {\em Solvability of a volume integral equation formulation for
  anisotropic elastodynamic scattering}, The Journal of Integral Equations and
  Applications, 28 (2016), pp.~169--203.

\bibitem{bonnet17}
\leavevmode\vrule height 2pt depth -1.6pt width 23pt, {\em A modified volume
  integral equation for anisotropic elastic or conducting inhomogeneities:
  unconditional solvability by neumann series}, J. Integral Equations Appl., 29
  (2017), pp.~271--295.

\bibitem{bonnet05}
{\sc M.~Bonnet and A.~Constantinescu}, {\em Inverse problems in elasticity},
  Inverse problems, 21 (2005), p.~R1.

\bibitem{ChaGinMin18}
{\sc R.~Chapko, D.~Gintides, and L.~Mindrinos}, {\em The inverse scattering
  problem by an elastic inclusion}, Advances in Computational Mathematics, 44
  (2018), pp.~453--476.

\bibitem{ChaKreMon00}
{\sc R.~Chapko, R.~Kress, and L.~M\"onch}, {\em On the numerical solution of a
  hypersingular integral equation for elastic scattering from a planar crack},
  IMA J. Numer. Anal., 20 (2000), pp.~601--619.

\bibitem{ClenCur60}
{\sc C.~W. Clenshaw and A.~R. Curtis}, {\em A method for numerical integration
  on an automatic computer}, Numerische Mathematik, 2 (1960), pp.~197--205.

\bibitem{ColKre13}
{\sc D.~Colton and R.~Kress}, {\em Inverse acoustic and electromagnetic
  scattering theory}, vol.~93 of Applied Mathematical Sciences, Springer, New
  York, 3~ed., 2013.

\bibitem{corna}
{\sc R.~Cornaggia}, {\em Development and use of higher-order asymptotics to
  solve inverse scattering problems}, PhD thesis, University of Minnesota,
  2016.

\bibitem{Doy00}
{\sc M.~M. Doyley, P.~M. Meaney, and J.~C. Bamber}, {\em Evaluation of an
  iterative reconstruction method for quantitative elastography}, Physics in
  Medicine \& Biology, 45 (2000), p.~1521.

\bibitem{eshelby1957}
{\sc J.~D. Eshelby}, {\em The determination of the elastic field of an
  ellipsoidal inclusion, and related problems}, Proceedings of the royal
  society of London. Series A. Mathematical and physical sciences, 241 (1957),
  pp.~376--396.

\bibitem{eshelby1961}
\leavevmode\vrule height 2pt depth -1.6pt width 23pt, {\em Elastic inclusions
  and inhomogeneities}, in Progress in Solid Mechanics, I.~N. Sneddon and
  R.~Hill, eds., vol.~2, North Holland, Amsterdam, 1961, pp.~89--140.

\bibitem{Feh06}
{\sc J.~Fehrenbach, M.~Masmoudi, R.~Souchon, and P.~Trompette}, {\em Detection
  of small inclusions by elastography}, Inverse problems, 22 (2006), p.~1055.

\bibitem{Gen72a}
{\sc W.~M. Gentleman}, {\em Implementing clenshaw-curtis quadrature, i
  methodology and experience}, Communications of the ACM, 15 (1972),
  pp.~337--342.

\bibitem{GinKir15}
{\sc D.~Gintides and K.~Kiriaki}, {\em Solvability of the integrodifferential
  equation of eshelby's equivalent inclusion method}, Quarterly Journal of
  Mechanics and Applied Mathematics, 68 (2015), pp.~85--96.

\bibitem{Hubmer}
{\sc S.~Hubmer, E.~Sherina, A.~Neubauer, and O.~Scherzer}, {\em Lam{\'e}
  parameter estimation from static displacement field measurements in the
  framework of nonlinear inverse problems}, SIAM Journal on Imaging Sciences,
  11 (2018), pp.~1268--1293.

\bibitem{ikehata}
{\sc M.~Ikehata}, {\em Inversion formulas for the linearized problem for an
  inverse boundary value problem in elastic prospection}, SIAM Journal on
  Applied Mathematics, 50 (1990), pp.~1635--1644.

\bibitem{Kirsch16}
{\sc A.~Kirsch and A.~Rieder}, {\em Inverse problems for abstract evolution
  equations with applications in electrodynamics and elasticity}, Inverse
  Problems, 32 (2016), p.~085001.

\bibitem{Kup79}
{\sc V.~Kupradze, T.~G. Gegelia, M.~O. Basheleshvili, and T.~V. Burchuladze},
  {\em Three-Dimensional problems of the mathematical theory of elasticity and
  thermoelasticity}, North-Holland Publishing Co., New York, 1979.

\bibitem{armin15}
{\sc A.~Lechleiter and J.~W. Schlasche}, {\em Identifying lam{\'e} parameters
  from time-dependent elastic wave measurements}, Inverse Problems in Science
  and Engineering, 25 (2017), pp.~2--26.

\bibitem{li2008}
{\sc S.~Li and G.~Wang}, {\em Introduction to micromechanics and
  nanomechanics}, World Scientific, 2008.

\bibitem{MarKouCha24}
{\sc V.~Markaki, D.~Kourounis, and A.~Charalambopoulos}, {\em On the
  identification of lam{\'e} parameters in linear isotropic elasticity via a
  weighted self-guided tv-regularization method}, Journal of Inverse and
  Ill-posed Problems, 32 (2024), pp.~213--231.

\bibitem{mikhlin87}
{\sc S.~G. Mikhlin and S.~Pr{\"o}ssdorf}, {\em Singular integral operators},
  vol.~68, Springer Science \& Business Media, 1987.

\bibitem{mura87}
{\sc T.~Mura}, {\em Micromechanics of defects in solids}, Martinus Nijhoff
  Publishers, The Netherlands, 2~ed., 1987.

\bibitem{nakamura93}
{\sc G.~Nakamura and G.~Uhlmann}, {\em Identification of lam{\'e} parameters by
  boundary measurements}, American Journal of Mathematics,  (1993),
  pp.~1161--1187.

\bibitem{nakamura94}
\leavevmode\vrule height 2pt depth -1.6pt width 23pt, {\em Global uniqueness
  for an inverse boundary problem arising in elasticity}, Inventiones
  mathematicae, 118 (1994), pp.~457--474.

\bibitem{nakamura95}
\leavevmode\vrule height 2pt depth -1.6pt width 23pt, {\em Inverse problems at
  the boundary for an elastic medium}, SIAM journal on mathematical analysis,
  26 (1995), pp.~263--279.

\bibitem{niwa1986}
{\sc Y.~Niwa, S.~Hirose, and M.~Kitahara}, {\em Elastodynamic analysis of
  inhomogeneous anisotropic bodies}, International journal of solids and
  structures, 22 (1986), pp.~1541--1555.

\bibitem{vladi}
{\sc V.~S. Vladimirov}, {\em Equations of mathematical physics}, Mir, 1984.

\end{thebibliography}

\end{document}